# COMPLEXITY REGULARIZATION VIA LOCALIZED RANDOM PENALTIES

By Gábor Lugosi[1] and Marten Wegkamp

*Pompeu Fabra University and Florida State University*

In this article, model selection via penalized empirical loss minimization in nonparametric classification problems is studied. Data-dependent penalties are constructed, which are based on estimates of the complexity of a small subclass of each model class, containing only those functions with small empirical loss. The penalties are novel since those considered in the literature are typically based on the entire model class. Oracle inequalities using these penalties are established, and the advantage of the new penalties over those based on the complexity of the whole model class is demonstrated.

**1. Introduction.** In this article, we propose a new complexity-penalized model selection method based on data-dependent penalties. We consider the binary classification problem where, given a random observation $X \in \mathbb{R}^d$, one has to predict $Y \in \{0, 1\}$. A *classifier* or *classification rule* is a function $f : \mathbb{R}^d \to \{0, 1\}$, with *loss*

$$L(f) \stackrel{\text{def}}{=} \mathbb{P}\{f(X) \neq Y\}.$$

A sample $\mathcal{D}_n = (X_1, Y_1), \ldots, (X_n, Y_n)$ of $n$ independent, identically distributed (i.i.d.) pairs is available. Each pair $(X_i, Y_i)$ has the same distribution as $(X, Y)$ and $\mathcal{D}_n$ is independent of $(X, Y)$. The statistician's task is to select a classification rule $f_n$ based on the data $\mathcal{D}_n$ such that the *probability of error*

$$L(f_n) = \mathbb{P}\{f_n(X) \neq Y | \mathcal{D}_n\}$$

is small. The Bayes classifier

$$f^*(x) \stackrel{\text{def}}{=} \mathbb{I}\{\mathbb{P}[Y = 1 | X = x] \geq \mathbb{P}[Y = 0 | X = x]\}$$

Received April 2002; revised July 2003.
[1]Supported by DGI Grant BMF2000-0807.

(where $\mathbb{I}$ denotes the indicator function) is the optimal rule as

$$L^* \stackrel{\text{def}}{=} \inf_{f:\mathbb{R}^d \to \{0,1\}} L(f) = L(f^*),$$

but both $f^*$ and $L^*$ are unknown to the statistician. In this article, we study classifiers $f:\mathbb{R}^d \to \{0,1\}$ which minimize the *empirical loss*

$$\widehat{L}(f) = \frac{1}{n}\sum_{i=1}^{n}\mathbb{I}\{f(X_i) \neq Y_i\}$$

over a class of rules $\mathcal{F}$. For any $\bar{f} \in \mathcal{F}$ minimizing the empirical probability of error, we have

$$\mathbb{E}L(\bar{f}) - L^* = \mathbb{E}\widehat{L}(\bar{f}) - L^* + \mathbb{E}(L - \widehat{L})(\bar{f})$$
$$= \mathbb{E}\inf_{f\in\mathcal{F}}\widehat{L}(f) - L^* + \mathbb{E}(L - \widehat{L})(\bar{f})$$
$$\leq \inf_{f\in\mathcal{F}}\mathbb{E}\widehat{L}(f) - L^* + \mathbb{E}(L - \widehat{L})(\bar{f})$$
$$= \inf_{f\in\mathcal{F}}L(f) - L^* + \mathbb{E}(L - \widehat{L})(\bar{f}).$$

Clearly, the *approximation error*

$$\inf_{f\in\mathcal{F}}L(f) - L^*$$

is decreasing as $\mathcal{F}$ becomes richer. However, the more complex $\mathcal{F}$, the more difficult the statistical problem becomes: the *estimation error*

$$\mathbb{E}(L - \widehat{L})(\bar{f})$$

increases with the complexity of $\mathcal{F}$. In many approaches to the problem described above, one fixes in advance a sequence of model classes $\mathcal{F}_1, \mathcal{F}_2, \ldots,$ whose union is $\mathcal{F}$. Denote by $\hat{f}_k$ a function in $\mathcal{F}_k$ having minimal empirical loss and by $L_k^* = \inf_{f\in\mathcal{F}_k} L(f)$ the minimal loss in class $\mathcal{F}_k$. The problem of penalized model selection is to find a possibly data-dependent penalty $\widehat{C}_k$, assigned to each class $\mathcal{F}_k$, such that minimizing the penalized empirical loss

$$\widehat{L}(f) + \widehat{C}_k, \qquad f \in \mathcal{F}_k,\ k=1,2,\ldots,$$

leads to a prediction rule

$$\hat{f} \stackrel{\text{def}}{=} \hat{f}_{\hat{k}}, \qquad \text{where } \hat{k} \stackrel{\text{def}}{=} \arg\min_{k\geq 1}(\widehat{L}(\hat{f}_k) + \widehat{C}_k),$$

with smallest possible loss.

The main idea is that since $\hat{f}_k$ minimizes $\widehat{L}(f)$ over $f \in \mathcal{F}_k$, we find, by the argument described above, that

$$\mathbb{E}L(\hat{f}_k) - L^* \leq L_k^* - L^* + \mathbb{E}(L - \widehat{L})(\hat{f}_k).$$

Our goal is to find the class $\mathcal{F}_k$ such that $L(\hat{f}_k)$ is as small as possible. To this end, a good balance has to be found between the approximation and estimation errors. The approximation error is unknown to us, but the estimation error may be estimated. The key to complexity-regularized model selection is that a tight bound for the estimation error is a good penalty $\widehat{C}_k$. More precisely, we show in Lemma 2.1 that if, for some constant $\gamma > 0$,

$$\mathbb{P}\{\widehat{C}_k \leq (L - \widehat{L})(\hat{f}_k)\} \leq \frac{\gamma}{n^2 k^2},$$

then the oracle inequality

$$\mathbb{E} L(\hat{f}) - L^* \leq \inf_k (L_k^* - L^* + \mathbb{E}\widehat{C}_k) + 2\gamma n^{-2}$$

holds, and also a similar bound,

$$L(\hat{f}) - L^* \leq \inf_k (L_k^* - L^* + 2\widehat{C}_k),$$

holds with probability greater than $1 - 4\gamma n^{-2}$. This simple result shows that the penalty should be, with large probability, an upper bound on the estimation error, and to guarantee good performance the bound should be as tight as possible.

Originally, distribution-free bounds, based on uniform-deviation inequalities, were proposed as penalties. For example, the structural risk minimization method of Vapnik and Chervonenkis [27] uses penalties of the form

$$\widehat{C}_k = \gamma \sqrt{\frac{\log \mathbb{S}_k(2n) + \log k}{n}},$$

where $\gamma$ is a constant and $\mathbb{S}_k(2n)$ is the $2n$-maximal shatter coefficient of the class

$$\mathcal{A}_k = \{\{x : f(x) = 1\}, \ f \in \mathcal{F}_k\},$$

that is,

$$\begin{aligned}
\mathbb{S}_k(2n) &= \max_{x_1,\ldots,x_{2n}} |\{\{x_1,\ldots,x_{2n}\} \cap A, \ A \in \mathcal{A}_k\}| \\
&= \max_{x_1,\ldots,x_{2n}} |\{(f(x_1),\ldots,f(x_{2n})), \ f \in \mathcal{F}_k\}|;
\end{aligned}$$
(1.1)

see, for example, [9, 26]. The fact that this type of penalty works follows from the Vapnik–Chervonenkis inequality. Such distribution-free bounds are attractive because of their simplicity, but precisely because of their distribution-free nature they are necessarily loose in many cases.

Recently, various attempts have been made to define the penalties in a data-dependent way to achieve this goal; see, for example, [2, 11, 13, 15, 17, 19, 22].

For example, in [2] and [11] random complexity penalties based on *Rademacher averages* were proposed and investigated. Rademacher averages are defined as

$$\widehat{R}_{\mathcal{F}_k} = \mathbb{E}\left[\sup_{f \in \mathcal{F}_k} \frac{1}{n}\sum_{i=1}^{n} \sigma_i \mathbb{I}\{f(X_i) \neq Y_i\} \bigg| \mathcal{D}_n \right],$$

where $\sigma_1, \ldots, \sigma_n$ are i.i.d. symmetric $\{-1, 1\}$-valued random variables independent of $\mathcal{D}_n$. The reason why this penalty was introduced is based on the fact that

$$\mathbb{E} \sup_{f \in \mathcal{F}_k} (L - \widehat{L})(f) \asymp \mathbb{E}\widehat{R}_{\mathcal{F}_k}$$

(see, e.g., [25]), and since $\widehat{R}_{\mathcal{F}_k}$ can be shown to be sharply concentrated around its mean. In fact, concentration inequalities have been a key tool in the analysis of data-based penalties (see [19]) and this paper relies heavily on some recent concentration results.

The model selection method based on Rademacher complexities satisfies an oracle inequality of the rough form

$$(1.2) \qquad \mathbb{E}L(\hat{f}) - L^* \leq \inf_k \left[ L_k^* - L^* + \gamma_1 \mathbb{E}\widehat{R}_{\mathcal{F}_k} + \gamma_2 \sqrt{\frac{\log k}{n}} \right]$$

(see [2] and [11]) for values of the constants $\gamma_1, \gamma_2 > 0$. The advantage of this bound over the one obtained by the distribution-free penalties mentioned above may perhaps be better understood if we further bound

$$\mathbb{E}\widehat{R}_{\mathcal{F}_k} \leq \sqrt{\frac{\mathbb{E}\log 2\mathbb{S}_k(X_1^n)}{2n}},$$

where

$$(1.3) \quad \begin{aligned} \mathbb{S}_k(X_1^n) &= |\{\{X_1, \ldots, X_n\} \cap A : A = \{x : f(x) = 1\}, \ f \in \mathcal{F}_k\}| \\ &= |\{(f(X_1), \ldots, f(X_n)), \ f \in \mathcal{F}_k\}|, \end{aligned}$$

is the *random shatter coefficient* of the class $\widehat{\mathcal{F}}_k$, which obviously never exceeds the worst-case shatter coefficient $\mathbb{S}_k(n)$ and may be significantly smaller for certain distributions.

However, this improved penalty is still not completely satisfactory. To see this, recall that by a classical result of Vapnik and Chervonenkis, for any index $k$,

$$(1.4) \qquad \mathbb{E}L(\hat{f}_k) - L_k^* \leq c\left( \sqrt{\frac{L_k^* \cdot \mathbb{E}\log \mathbb{S}_k(X_1^n)}{n}} + \frac{\mathbb{E}\log \mathbb{S}_k(X_1^n)}{n} \right),$$

which is much smaller than the corresponding expected Rademacher average if $L_k^*$ is small. (For explicit constants we refer to Theorem 1.14 in [16].) Since in typical classification problems the minimal error $L_k^*$ in class $\mathcal{F}_k$ is often very small for some $k$, it is important to find penalties which allow derivation of oracle inequalities with the appropriate dependence on $L_k^*$. In particular, a desirable goal would be to develop classifiers $\hat{f}$ for which an oracle inequality resembling

$$\mathbb{E}L(\hat{f}) - L^* \leq \inf_k \left\{ L_k^* - L^* + \gamma_1 \sqrt{\frac{L_k^* \cdot \mathbb{E}\log \mathbb{S}_k(X_1^n)}{n}} + \gamma_2 \frac{\mathbb{E}\log \mathbb{S}_k(X_1^n)}{n} \right\}$$

holds for all distributions. The main results of this article (Theorems 4.1 and 4.2) show that estimates of the desired property are indeed possible to construct in a conceptually simple way.

By the key Lemma 2.1, it suffices to find a data-dependent upper estimate of $(L - \widehat{L})(\hat{f}_k)$ which has the order of magnitude of the above upper bound. The difficulty is that $L_k^*$ and $\mathbb{E}\log \mathbb{S}_k(X_1^n)$ both depend on the underlying distribution.

The improvement is achieved by decreasing the penalties so that the supremum in the definition of the Rademacher average is not taken over the whole class $\mathcal{F}_k$ but rather over a small subclass $\widehat{\mathcal{F}}_k$ containing only functions which "look good" on the data. More precisely, define the random subclass $\widehat{\mathcal{F}}_k \subset \mathcal{F}_k$ by

$$\widehat{\mathcal{F}}_k = \{f \in \mathcal{F}_k : \widehat{L}(f) \leq \gamma_1 \widehat{L}(\hat{f}_k) + \gamma_2 n^{-1} \log \mathbb{S}_k(X_1^n) + \gamma_3 n^{-1} \log(nk)\}$$

for some nonnegative constants $\gamma_1, \gamma_2$ and $\gamma_3$.

Risk estimates based on localized Rademacher averages have been considered in several recent papers. The most closely related procedure is proposed by Koltchinskii and Panchenko [12], who, assuming $\inf_{f \in \mathcal{F}} L(f) = 0$, compute the Rademacher averages of subclasses of $\mathcal{F}$ with empirical loss less than $r$ for different values of $r$ obtained by a recursive procedure, and obtain bounds for the loss of the empirical risk minimizer in terms of the localized Rademacher averages obtained after a certain number of iterations. Our approach of bounding the loss is conceptually simpler: it suffices to compute the Rademacher complexities at only one scale which depends on the smallest empirical loss in the class and a term of a smaller order determined by the shatter coefficients of the whole class. Thus, we use "global" information to determine the scale of localization. Bartlett, Bousquet and Mendelson [3] also derive closely related generalization bounds based on localized Rademacher averages. In their approach the performance bounds also depend on Rademacher averages computed at different scales of localization, which are combined by the technique of peeling. For further recent related work, we also refer to [7, 8, 24].

The rest of the paper is organized as follows. Section 2 presents some basic inequalities on model selection, which generalizes some of the results in [2]. Section 3 proposes a simple but suboptimal penalty which already has some of the main features of the penalties presented in Section 4. It shows, in a transparent way, some of the underlying ideas of the main results. Section 4 introduces a new penalty based on the Rademacher average $\widehat{R}_{\widehat{\mathcal{F}}_k}$ and it is shown that the new estimate yields an improvement of the desired form.

**2. Preliminaries.** In this section we present two basic auxiliary lemmata on model selection. The first lemma is general in the sense that it does not depend on the particular choice of the penalty $\widehat{C}_k$. This result was mentioned in the Introduction and generalizes a result obtained by Bartlett, Boucheron and Lugosi [2]. We recall that the penalized estimator is defined by $\hat{f} \stackrel{\text{def}}{=} \hat{f}_{\hat{k}}$, with $\hat{k} \stackrel{\text{def}}{=} \arg\min_{k \geq 1}(\widehat{L}(\hat{f}_k) + \widehat{C}_k)$.

LEMMA 2.1. *Suppose that the random variables $\widehat{C}_1, \widehat{C}_2, \ldots$ are such that*

$$\mathbb{P}\{\widehat{C}_k \leq (L - \widehat{L})(\hat{f}_k)\} \leq \frac{\gamma}{n^2 k^2}$$

*for some $\gamma > 0$ and for all $k$. Then we have*

$$\mathbb{E}L(\hat{f}) - L^* \leq \inf_k [L_k^* - L^* + \mathbb{E}\widehat{C}_k] + \frac{2\gamma}{n^2}.$$

It is clear that we can always take $\widehat{C}_k \leq 1$.

PROOF OF LEMMA 2.1. Observe that

$$\mathbb{E}\sup_k\{(L - \widehat{L})(\hat{f}_k) - \widehat{C}_k\} \leq \mathbb{P}\left\{\sup_k[(L - \widehat{L})(\hat{f}_k) - \widehat{C}_k] \geq 0\right\}$$

$$(\text{since } \sup_k[(L - \widehat{L})(\hat{f}_k) - \widehat{C}_k] \leq 1)$$

$$\leq \sum_{k=1}^{\infty} \mathbb{P}\{(L - \widehat{L})(\hat{f}_k) - \widehat{C}_k \geq 0\}$$

$$(\text{by the union bound})$$

$$\leq \sum_{k=1}^{\infty} \frac{\gamma}{n^2 k^2}$$

$$(\text{by assumption})$$

$$\leq \frac{2\gamma}{n^2}.$$

Therefore, we may conclude that

$$\mathbb{E}L(\hat{f}) - L^* = \mathbb{E}[\widehat{L}(\hat{f}) - L^* + \widehat{C}_{\hat{k}}] + \mathbb{E}[(L - \widehat{L})(\hat{f}) - \widehat{C}_{\hat{k}}]$$

(where $\hat{k}$ is the selected model index, i.e., $\hat{f} = \hat{f}_{\hat{k}}$)

$$\leq \mathbb{E}\inf_k[\widehat{L}(\hat{f}_k) - L^* + \widehat{C}_k] + \mathbb{E}[(L - \widehat{L})(\hat{f}) - \widehat{C}_{\hat{k}}]$$

(by definition of $\hat{f}$)

$$\leq \mathbb{E}\inf_k\left[\inf_{f \in \mathcal{F}_k} \widehat{L}(f) - L^* + \widehat{C}_k\right] + \mathbb{E}\sup_k[(L - \widehat{L})(\hat{f}_k) - \widehat{C}_k]$$

(by definition of $\hat{f}_k$)

$$\leq \inf_k\left[\inf_{f \in \mathcal{F}_k} L(f) - L^* + \mathbb{E}\widehat{C}_k\right] + \mathbb{E}\sup_k[(L - \widehat{L})(\hat{f}_k) - \widehat{C}_k]$$

(interchange $\mathbb{E}$ and inf)

$$\leq \inf_k[L_k^* - L^* + \mathbb{E}\widehat{C}_k] + \frac{2\gamma}{n^2}$$

(by the preceding display)

and the proof is complete. $\square$

The preceding result is not entirely satisfactory for the following reason. Although it presents a useful bound, it is a bound for the *average* risk behavior of $\hat{f}$. However, the penalty is computed on the data at hand, and therefore the proposed criterion should have optimal performance for (almost) all possible sequences of the data. The following result presents a nonasymptotic oracle inequality which holds with large probability and an asymptotic almost-sure version.

LEMMA 2.2. *Assume that, for all $k, n \geq 1$,*

$$\mathbb{P}\{\widehat{C}_k \leq (L - \widehat{L})(\hat{f}_k)\} \leq \frac{\gamma}{n^2 k^2}$$

*and*

$$\mathbb{P}\{\widehat{C}_k \leq (\widehat{L} - L)(f_k^*)\} \leq \frac{\gamma}{n^2 k^2}.$$

*Then, for all $n \geq 1$ we have*

$$\mathbb{P}\left[L(\hat{f}) - L^* \geq \inf_k(L_k^* - L^* + 2\widehat{C}_k)\right] \leq \frac{4\gamma}{n^2}$$

*and the asymptotic almost-sure bound*

$$\mathbb{P}\left[\liminf_{n \to \infty}\left\{L(\hat{f}) - L^* \leq \inf_k(L_k^* - L^* + 2\widehat{C}_k)\right\}\right] = 1.$$

PROOF. Let $\hat{k}$ be the selected model index. Notice that

$$\begin{aligned}L(\hat{f}) &= \widehat{L}(\hat{f}) + \widehat{C}_{\hat{k}} + (L - \widehat{L})(\hat{f}) - \widehat{C}_{\hat{k}} \\ &\leq \inf_k [\widehat{L}(\hat{f}_k) + \widehat{C}_k] + \sup_k [(L - \widehat{L})(\hat{f}_k) - \widehat{C}_k] \\ &\leq \inf_k [\widehat{L}(f_k^*) + \widehat{C}_k] + \sup_k [(L - \widehat{L})(\hat{f}_k) - \widehat{C}_k] \\ &\leq \inf_k [L_k^* + 2\widehat{C}_k] + \sup_k [(\widehat{L} - L)(f_k^*) - \widehat{C}_k] + \sup_k [(L - \widehat{L})(\hat{f}_k) - \widehat{C}_k].\end{aligned}$$

By assumption, the last two terms on the right-hand side satisfy

$$\mathbb{P}\left[\sup_k [(\widehat{L} - L)(f_k^*) - \widehat{C}_k] + \sup_k [(L - \widehat{L})(\hat{f}_k) - \widehat{C}_k] \geq 0\right] \leq \sum_{k=1}^{\infty} \frac{2\gamma}{n^2 k^2} < \frac{4\gamma}{n^2},$$

proving the first inequality. The almost-sure statement is a direct consequence of the Borel–Cantelli lemma. □

**3. A simple version.** The purpose of this short section is to offer a simplified yet suggestive illustration of the ideas. As discussed in the Introduction, an ideal penalty would be a tight upper bound for the expression on the right-hand side of (1.4). Motivated by this bound, we propose the simple penalty

$$\widehat{C}_k = 2\sqrt{2\widehat{L}(\hat{f}_k) + 8\frac{\log \mathbb{S}_k(2n) + 2\log(nk)}{n}} \cdot \sqrt{\frac{\log \mathbb{S}_k(2n)}{n}} + 2\frac{\log(nk)}{n},$$

where $\mathbb{S}_k(2n)$ is the (worst-case) $2n$-shatter coefficient defined in (1.1). Thus, the minimal loss $L_k^*$ in class $\mathcal{F}_k$ is estimated by its natural empirical counterpart $\widehat{L}(\hat{f}_k) = \inf_{f \in \mathcal{F}_k} \widehat{L}(f)$ and the expected logarithmic shatter coefficient $\mathbb{E}\log \mathbb{S}_k(X_1^n)$ is estimated by the distribution-free upper bound $\log \mathbb{S}_k(2n)$. [This term may be bounded further by $V_k \log(2n + 1)$, where $V_k$ is the VC-dimension of the set $\mathcal{A}_k$.] The auxiliary terms $n^{-1}\log(nk)$ are necessary to derive the desired oracle inequalities. The next theorem shows that the proposed penalty indeed works.

THEOREM 3.1. *Consider the penalized empirical loss minimizer $\hat{f}$ with the data-based penalty $\widehat{C}_k$ defined above. Then, for every $n$ and for all distributions of $(X, Y)$,*

$$\mathbb{E}L(\hat{f}) - L^* \leq \inf_k (L_k^* - L^* + \mathbb{E}\widehat{C}_k) + \frac{16}{n^2}.$$

*In particular,*

$$\mathbb{E}L(\hat{f}) - L^* \leq \inf_k \left[ L_k^* - L^* + 4\sqrt{L_k^* + \frac{2}{n}\{\log \mathbb{S}_k(2n) + 2\log(nk)\}}\right.$$

RANDOM PENALTIES 9

$$\times \sqrt{\frac{\log \mathbb{S}_k(2n)}{n} + 2\frac{\log(nk)}{n}} \,\Bigg] + \frac{16}{n^2}.$$

The proof uses Lemma 2.1 and the following uniform deviation bound due to Vapnik and Chervonenkis [27]. (The slightly improved form used here is proved by Anthony and Shawe-Taylor [1].)

PROPOSITION 3.2. *Let $\mathbb{S}_k(X_1^{2n})$ be the random shatter coefficient of $\mathcal{A}_k$ based on i.i.d. observations $X_1,\ldots,X_{2n}$ defined in* (1.3). *For all $\varepsilon > 0$ and $n \geq 1$,*

(3.1) $$\mathbb{P}\Big\{\sup_{f \in \mathcal{F}} L(f) - 2\widehat{L}(f) \geq 2\varepsilon\Big\} \leq 4\mathbb{E}\mathbb{S}_k(X_1^{2n})\exp(-n\varepsilon/4)$$

*and*

(3.2) $$\mathbb{P}\Big\{\sup_{f \in \mathcal{F}} \widehat{L}(f) - 2L(f) \geq 2\varepsilon\Big\} \leq 4\mathbb{E}\mathbb{S}_k(X_1^{2n})\exp(-n\varepsilon/4).$$

PROOF. Observe that, for all $\varepsilon > 0$ and $n \geq 1$,

$$\Big\{\sup_{f \in \mathcal{F}} L(f) - 2\widehat{L}(f) \geq 2\varepsilon\Big\} \subseteq \Big\{\sup_{f \in \mathcal{F}} \frac{L(f) - \widehat{L}(f)}{\sqrt{L(f)}} \geq \sqrt{\varepsilon}\Big\},$$

and similarly,

$$\Big\{\sup_{f \in \mathcal{F}} \widehat{L}(f) - 2L(f) \geq 2\varepsilon\Big\} \subseteq \Big\{\sup_{f \in \mathcal{F}} \frac{\widehat{L}(f) - L(f)}{\sqrt{\widehat{L}(f)}} \geq \sqrt{\varepsilon}\Big\}.$$

The proposition follows by [1]. □

PROOF OF THEOREM 3.1. We start with the proof of the first inequality of Theorem 3.1. In view of Lemma 2.1, it suffices to show that

$$\mathbb{P}\{L(\hat{f}_k) - \widehat{L}(\hat{f}_k) \geq \widehat{C}_k\} \leq 8/(nk)^2.$$

Consequently, by (3.2),

$$\mathbb{P}\Big\{2\widehat{L}(\hat{f}_k) + 8\frac{\log \mathbb{S}_k(2n)}{n} + 16\frac{\log(nk)}{n} \leq L(\hat{f}_k)\Big\}$$

$$= \mathbb{P}\Big\{L(\hat{f}_k) - 2\widehat{L}(\hat{f}_k) \geq 8\frac{\log \mathbb{S}_k(2n)}{n} + 16\frac{\log(nk)}{n}\Big\}$$

$$\leq 4\mathbb{S}_k(2n)\exp\Big\{-\frac{n}{8}\Big(8\frac{\log \mathbb{S}_k(2n)}{n} + 16\frac{\log(nk)}{n}\Big)\Big\}$$

$$= \frac{4}{n^2 k^2},$$

so that
$$\mathbb{P}\{\widehat{C}_k \geq \widetilde{C}_k\} \geq 1 - 4/(nk)^2,$$
where
$$\widetilde{C}_k = 2\sqrt{L(\hat{f}_k)} \cdot \sqrt{\frac{\log \mathbb{S}_k(2n)}{n} + 2\frac{\log(nk)}{n}}.$$

Another application of inequality (3.2) yields

$$\mathbb{P}\{L(\hat{f}_k) - \widehat{L}(\hat{f}_k) \geq \widehat{C}_k\}$$
$$\leq \mathbb{P}\{L(\hat{f}_k) - \widehat{L}(\hat{f}_k) \geq \widetilde{C}_k\} + \frac{4}{(nk)^2}$$
$$\leq 4\mathbb{S}_k(2n) \exp\left\{-\frac{n}{4} \cdot 4\left(\frac{\log \mathbb{S}_k(2n)}{n} + 2\frac{\log(nk)}{n}\right)\right\} + \frac{4}{(nk)^2}$$
$$= \frac{8}{(nk)^2}.$$

Conclude via Lemma 2.1 that
$$\mathbb{E}L(\hat{f}) \leq \min_k (L_k + \mathbb{E}\widehat{C}_k) + \frac{16}{n^2}.$$

For the second inequality, deduce that for all $\delta > 0$,
$$\mathbb{E}\sqrt{\widehat{L}(\hat{f}_k) + \delta} \leq \sqrt{\mathbb{E}\widehat{L}(\hat{f}_k) + \delta} \leq \sqrt{\mathbb{E}\inf_{f \in \mathcal{F}_k} \widehat{L}(f) + \delta} \leq \sqrt{L_k^* + \delta},$$

by Jensen's inequality and the definition of $\hat{f}_k$. $\square$

The bound of Theorem 3.1 has the right dependence on $L_k^*$ as suggested by inequality (1.4) mentioned in the Introduction. In particular, if $L_k^*$ happens to equal zero for some class $\mathcal{F}_k$, then the upper bound has an improved rate of convergence. The disadvantage of the simple penalty defined above is that instead of the expected shatter coefficients, a distribution-free (and therefore suboptimal) upper bound appears for each class $\mathcal{F}_k$.

Recently, Boucheron, Lugosi and Massart [4] proved that $\log \mathbb{S}_k(X_1^n)$ concentrates sharply around its mean. For example, we have the following inequalities.

PROPOSITION 3.3. *For all $\varepsilon > 0$, $n \geq 1$,*
$$\mathbb{P}[\mathbb{E}\log \mathbb{S}_k(X_1^n) > 2\log \mathbb{S}_k(X_1^n) + 2\varepsilon] \leq e^{-\varepsilon},$$
$$\mathbb{P}[\log \mathbb{S}_k(X_1^n) > 2\mathbb{E}\log \mathbb{S}_k(X_1^n) + 2\varepsilon] \leq e^{-\varepsilon}.$$

*Moreover, for each $n \geq 1$,*
$$\mathbb{E}\log \mathbb{S}_k(X_1^n) \leq \log \mathbb{E}\mathbb{S}_k(X_1^n) \leq \frac{1}{\ln 2}\mathbb{E}\log \mathbb{S}_k(X_1^n) \leq 2\mathbb{E}\log \mathbb{S}_k(X_1^n).$$

This proposition implies that the expected random log shatter coefficients $\mathbb{E}\log\mathbb{S}_k(X_1^n)$ of $\mathcal{F}_k$ may be replaced by a constant times $\log\mathbb{S}_k(X_1^n)$ and vice versa. Hence we may replace the distribution-free bounds $\log\mathbb{S}_k(2n)$ by empirical estimates $\log\mathbb{S}_k(X_1^n)$, at the price of slightly worse constants. The main oracle inequalities in Section 4 are accompanied by asymptotic almost-sure versions of bounds for the expected value. Such bounds are easy to obtain as well, simply by invoking Lemma 2.2 instead of Lemma 2.1. The details are omitted here.

**4. Rademacher penalties.** The main results of the paper are presented in this section. Assign to each model class $\mathcal{F}_k$,

$$(4.1) \qquad \hat{u}_k = 16\frac{4\log\mathbb{S}_k(X_1^n) + 9\log(nk)}{n},$$

with $\mathbb{S}_k(X_1^n)$ defined in (1.3), and the class

$$(4.2) \qquad \widehat{\mathcal{F}}_k = \{f \in \mathcal{F}_k : \widehat{L}(f) \leq 16\widehat{L}(\hat{f}_k) + 15\hat{u}_k\}.$$

Observe that the class $\widehat{\mathcal{F}}_k$ contains only those classifiers whose empirical loss is not much larger than that of the empirical minimizer. Note that the constant 16 has no special role; it has been chosen by convenience. Any constant larger than 1 would lead to similar results, at the price of modifying other constants. The term $\hat{u}_k$ depends on the shatter coefficient of the whole class $\mathcal{F}_k$ but it is typically small compared to $\widehat{L}(\hat{f}_k)$.

The penalty is calculated in terms of the Rademacher average of this smaller class. More precisely, define the complexity estimate by

$$(4.3) \quad \widehat{C}_k = (8\widehat{R}_{\widehat{\mathcal{F}}_k} + 20n^{-1}\log(nk) + 2\sqrt{n^{-1}\log(nk)} \cdot \sqrt{8\widehat{L}(\hat{f}_k) + 7\hat{u}_k}) \wedge 1.$$

Again, not too much attention should be paid to the values of the constants involved. We favored simple readable proofs over optimal constants. Note that, through $\mathbb{S}_k(X_1^n)$, the penalty also depends on the random shatter coefficient of the whole class $\mathcal{F}_k$. However, the term involving the shatter coefficient of the entire class $\mathcal{F}_k$,

$$n^{-1}\sqrt{\log(nk) \cdot \log\mathbb{S}_k(X_1^n)},$$

is typically much smaller (by a factor $n^{-1/2}$) than the Rademacher average of the whole class $\mathcal{F}_k$. [For instance, see (4.8) and Proposition 4.6.]

We have the following performance bound for the expected loss of the minimizer $\hat{f}$ of the penalized empirical loss $\widehat{L}(\hat{f}_k) + \widehat{C}_k$.

THEOREM 4.1. *For every $n$,*

$$\mathbb{E}L(\hat{f}) - L^* \leq \inf_k(L_k^* - L^* + \mathbb{E}\widehat{C}_k) + \frac{22}{n^2}.$$

*In addition, with probability greater than $1 - 44/n^2$,*

$$L(\hat{f}) - L^* \leq \inf_k (L_k^* - L^* + 2\widehat{C}_k),$$

*and also*

$$\mathbb{P}\left[\liminf_{n \to \infty}\left\{L(\hat{f}) - L^* \leq \inf_k (L_k^* - L^* + 2\widehat{C}_k)\right\}\right] = 1.$$

The next theorem is here to point out that the bound above is indeed a significant improvement over bounds of the type (1.2), and that the dependence on the minimal loss $L_k^*$ and the random shatter coefficient has the form suggested by (1.4). For this purpose, we introduce

(4.4) $$\bar{u}_k = 16 \frac{8 \mathbb{E} \log \mathbb{S}_k(X_1^n) + 17 \log(nk)}{n}$$

and the class

$$\overline{\mathcal{F}}_k = \{f \in \mathcal{F}_k : L(f) \leq 64 L_k^* + 63 \bar{u}_k\}.$$

We also set

$$\varepsilon_k = 2n^{-1} \log(nk).$$

THEOREM 4.2. *The following oracle inequality holds:*

$$\mathbb{E} L(\hat{f}) - L^* \leq \min_{k \geq 1}[L_k^* - L^* + 8 \mathbb{E} \widehat{R}_{\overline{\mathcal{F}}_k} + 15\varepsilon_k + 16\sqrt{L_k^* + \bar{u}_k} \cdot \sqrt{2\varepsilon_k}] + 22n^{-2}.$$

*In particular, there exist universal constants $\gamma_1$ and $\gamma_2$ such that*

$$\mathbb{E} L(\hat{f}) - L^* \leq \inf_k \left\{ L_k^* - L^* + \gamma_1 \sqrt{\frac{L_k^* \cdot (\mathbb{E} \log \mathbb{S}_k(X_1^n) \vee \log(nk))}{n}} \right.$$
$$\left. + \gamma_2 \frac{\mathbb{E} \log \mathbb{S}_k(X_1^n) \vee \log(nk)}{n} \right\}.$$

This oracle inequality has the desired form outlined in the Introduction and improves upon the results of [2] and [13]. For example, in the special case when $L_k^* = 0$ for $k \geq k_0$, we obtain, for some numerical constants $c_1$ and $c_2$,

$$\mathbb{E} L(\hat{f}) \leq \min_{k \geq k_0} c_1 \frac{\mathbb{E} \log \mathbb{S}_k(X_1^n) \vee \log(nk)}{n} + \frac{c_2}{n^2},$$

which is of a different order of magnitude from the penalties considered by [2] and [13]. Theorem 4.2 is only stated for the expected loss but an inequality which holds with "large" probability may be obtained just as in Theorem 4.1.

*Proofs of Theorems* 4.1 *and* 4.2. First, recall the definitions of $\hat{u}_k$ and $\bar{u}_k$ in (4.1) and (4.4), respectively, and in addition define

$$u_k = 8 \frac{2 \log \mathbb{E}\mathbb{S}_k(X_1^n) + 2 \log(nk)}{n}$$

and the event

$$B_k \stackrel{\text{def}}{=} \{u_k \leq \hat{u}_k \leq \bar{u}_k\}.$$

Observe Proposition 3.3 yields that, with probability at least $1 - 1/(nk)^2$,

$$\begin{aligned}
u_k &= \frac{16}{n}\{\log \mathbb{E}\mathbb{S}_k(X_1^n) + \log(nk)\} \\
&\leq \frac{16}{n}\{2\mathbb{E}\log \mathbb{S}_k(X_1^n) + \log(nk)\} \\
&\leq \frac{16}{n}\{2[2\log \mathbb{S}_k(X_1^n) + 4\log(nk)] + \log(nk)\} \\
&= \hat{u}_k \\
&\leq \frac{16}{n}\{4[2\mathbb{E}\log \mathbb{S}_k(X_1^n) + 4\log(nk)] + 9\log(nk)\} \\
&= \frac{16}{n}\{8\mathbb{E}\log \mathbb{S}_k(X_1^n) + 17\log(nk)\} \\
&= \bar{u}_k
\end{aligned}$$

and therefore

(4.5) $$\mathbb{P}B_k^c \leq (nk)^{-2}.$$

Finally, we introduce the event

$$A_k = \left\{\sup_{f \in \mathcal{F}_k} L(f) - 2\widehat{L}(f) \leq u_k\right\} \cap \left\{\sup_{f \in \mathcal{F}_k} \widehat{L}(f) - 2L(f) \leq u_k\right\}$$

and the class

$$\mathcal{F}_k^* = \{f \in \mathcal{F}_k \colon L(f) \leq 4L_k^* + 3u_k\}.$$

The following intermediate result will be useful in the proofs of both theorems.

LEMMA 4.3. *We have*

(4.6) $$\mathbb{P}\{A_k \cap B_k\} \geq 1 - \frac{9}{(nk)^2},$$

*and on the set* $A_k \cap B_k$ *the following hold:*

(i) $\hat{f}_k \in \mathcal{F}_k^*$.

(ii) $\mathcal{F}_k^* \subseteq \widehat{\mathcal{F}}_k$, and in particular, $\widehat{R}_{\mathcal{F}_k^*} \leq \widehat{R}_{\widehat{\mathcal{F}}_k}$.
(iii) $L_k^* \leq 2\widehat{L}(\hat{f}_k) + u_k$.

PROOF. To begin with, notice that

$$\mathbb{E}\mathbb{S}_k(X_1^{2n}) \leq \mathbb{E}\mathbb{S}_k(X_1^n)\mathbb{S}_k(X_{n+1}^{2n}) = \mathbb{E}^2\mathbb{S}_k(X_1^n)$$

by the definition of the shatter coefficient and by the independence of the $X_i$. Thus, by Proposition 3.2,

$$\mathbb{P}A_k^c \leq 8\mathbb{E}\mathbb{S}_k(X_1^{2n})\exp\left(-\frac{nu_k}{8}\right) \leq \frac{8}{n^2k^2}.$$

This bound and (4.5) imply assertion (4.6). To prove claim (i), observe that on $A_k$,

$$\begin{aligned}
L(\hat{f}_k) &\leq 2\widehat{L}(\hat{f}_k) + u_k & &\text{(by definition of } A_k\text{)} \\
&\leq 2\widehat{L}(f_k^*) + u_k & &\text{(by definition of } \hat{f}_k\text{)} \\
&\leq 2(2L_k^* + u_k) + u_k & &\text{(by definition of } A_k\text{)} \\
&= 4L_k^* + 3u_k.
\end{aligned}$$

For claim (ii), notice that, for any $f \in \mathcal{F}_k^*$,

$$\begin{aligned}
\widehat{L}(f) &\leq 2L(f) + u_k & &\text{(by definition of } A_k\text{)} \\
&\leq 2[4L_k^* + 3u_k] + u_k & &\text{(by definition of } \mathcal{F}_k^*\text{)} \\
&= 8L_k^* + 7u_k \\
&\leq 8L(\hat{f}_k) + 7u_k & &\text{(by definition of } L_k^*\text{)} \\
&\leq 16\widehat{L}(\hat{f}_k) + 15u_k & &\text{(by definition of } A_k\text{)} \\
&\leq 16\widehat{L}(\hat{f}_k) + 15\hat{u}_k & &\text{(by definition of } B_k\text{)}.
\end{aligned}$$

Claim (ii) now follows. Claim (iii) is immediate from the definition of $A_k$ since both $\hat{f}_k$ and $f_k^*$ belong to $\mathcal{F}_k$. □

Next we link the Rademacher average $\widehat{R}_{\mathcal{F}_k^*}$ to $\mathbb{E}\sup_{f \in \mathcal{F}_k^*} |L(f) - \widehat{L}(f)|$. By a classical symmetrization device (cf. [10] or [25]),

$$(4.7) \qquad \mathbb{E}\sup_{f \in \mathcal{F}_k^*} |\widehat{L}(f) - L(f)| \leq 2\mathbb{E}\widehat{R}_{\mathcal{F}_k^*}.$$

Also, $\widehat{R}_{\mathcal{F}_k^*}$ is known to concentrate sharply around its mean. For example, we have, by results of [4, 5], the following bounds.

PROPOSITION 4.4. *For all $\epsilon > 0$, $n \geq 1$,*

$$\mathbb{P}[\widehat{R}_{\mathcal{F}_k} \geq 2\mathbb{E}\widehat{R}_{\mathcal{F}_k} + \epsilon] \leq e^{-6n\epsilon/5} \quad and \quad \mathbb{P}[\widehat{R}_{\mathcal{F}_k} \leq \tfrac{1}{2}\mathbb{E}\widehat{R}_{\mathcal{F}_k} - \epsilon] \leq e^{-n\epsilon}.$$

PROOF. Define $Z \stackrel{\text{def}}{=} n\widehat{R}_{\mathcal{F}_k}$. Then it follows from [4] that

$$\log \mathbb{E}\exp(\lambda(Z - \mathbb{E}Z)) \leq \mathbb{E}Z(e^\lambda - 1 - \lambda),$$

which implies further that, for $0 \leq \lambda < 3$,

$$\log \mathbb{E}\exp(\lambda(Z - \mathbb{E}Z)) \leq \frac{\lambda \mathbb{E}Z}{2(1 - \lambda/3)}.$$

After an application of Markov's inequality, we find

$$\mathbb{P}[Z \geq \mathbb{E}Z + \sqrt{2\mathbb{E}Zx} + x/3] \leq e^{-x}.$$

We obtain the desired upper-tail bound by inserting $Z = n\widehat{R}_{\mathcal{F}_k}$ in the preceding display and invoking the inequality $2\sqrt{xy} \leq x + y$. The bound for the lower tail follows from the inequality

$$\mathbb{P}[Z \leq \mathbb{E}Z - \sqrt{2x\mathbb{E}Z}] \leq e^{-x}$$

(see [4]) and since $x + \tfrac{1}{2}y \geq \sqrt{2xy}$. □

Finally, we make key use of the following concentration inequality for the supremum of an empirical process, recently established by Talagrand [23]; see also [14, 19, 21]. The best-known constants reported here have been obtained by Bousquet [6].

PROPOSITION 4.5. *Set $\Sigma_{\mathcal{F}_k^*} = \sup_{f \in \mathcal{F}_k^*} L(f)(1 - L(f))$. For all $\epsilon > 0$, $n \geq 1$,*

$$\mathbb{P}\left[\sup_{f \in \mathcal{F}_k^*} |\widehat{L}(f) - L(f)| \geq 2\mathbb{E}\sup_{f \in \mathcal{F}_k^*} |\widehat{L}(f) - L(f)| + \Sigma_{\mathcal{F}_k^*}\sqrt{2\epsilon} + \frac{4\epsilon}{3}\right] \leq e^{-n\epsilon}.$$

We are now ready to prove Theorems 4.1 and 4.2.

PROOF OF THEOREM 4.1. Deduce, using (i), (ii) and (iii) of Lemma 4.3, the following string of inequalities:

$$\mathbb{P}[\{L(\hat{f}_k) \geq \widehat{L}(\hat{f}_k) + \widehat{C}_k\} \cap A_k \cap B_k]$$
$$= \mathbb{P}[\{L(\hat{f}_k) \geq \widehat{L}(\hat{f}_k) + 8\widehat{R}_{\widehat{\mathcal{F}}_k}$$
$$\qquad + 10\varepsilon_k + \sqrt{8\widehat{L}(\hat{f}_k) + 7\hat{u}_k}\sqrt{2\varepsilon_k}\} \cap A_k \cap B_k]$$
$$\leq \mathbb{P}[\{\exists f \in \mathcal{F}_k^* : L(f) \geq \widehat{L}(f) + 8\widehat{R}_{\widehat{\mathcal{F}}_k}$$

$$+ 10\varepsilon_k + \sqrt{8\widehat{L}(\hat{f}_k) + 7\widehat{u}_k}\sqrt{2\varepsilon_k}\} \cap A_k \cap B_k]$$

[by property (i)]

$$\leq \mathbb{P}[\{\exists f \in \mathcal{F}_k^* : L(f) \geq \widehat{L}(f) + 8\widehat{R}_{\mathcal{F}_k^*}$$
$$+ 10\varepsilon_k + \sqrt{8\widehat{L}(\hat{f}_k) + 7u_k}\sqrt{2\varepsilon_k}\} \cap A_k \cap B_k]$$

[by property (ii) and definition of $B_k$]

$$\leq \mathbb{P}[\{\exists f \in \mathcal{F}_k^* : L(f) \geq \widehat{L}(f) + 8\widehat{R}_{\mathcal{F}_k^*}$$
$$+ 10\varepsilon_k + \sqrt{4L_k^* + 3u_k}\sqrt{2\varepsilon_k}\} \cap A_k \cap B_k]$$

[by property (iii)]

$$\leq \mathbb{P}\left\{\sup_{f \in \mathcal{F}_k^*} |L(f) - \widehat{L}(f)| \geq 8\widehat{R}_{\mathcal{F}_k^*} + 10\varepsilon_k + \Sigma_{\mathcal{F}_k^*}\sqrt{2\varepsilon_k}\right\},$$

where the last inequality follows from

$$\Sigma_{\mathcal{F}_k^*}^2 = \sup_{f \in \mathcal{F}_k^*} \mathrm{Var}(\mathbb{I}\{f(X) \neq Y\}) \leq \sup_{f \in \mathcal{F}_k^*} L(f) \leq 4L_k^* + 3u_k.$$

Invoke (4.7), (4.6) and Propositions 4.4 and 4.5 to conclude that

$$\mathbb{P}\{L(\hat{f}_k) \geq \widehat{L}(\hat{f}_k) + \widehat{C}_k\}$$
$$\leq \mathbb{P}\left\{\sup_{f \in \mathcal{F}_k^*} |L(f) - \widehat{L}(f)| \geq 8\widehat{R}_{\mathcal{F}_k^*} + 10\varepsilon_k + \Sigma_{\mathcal{F}_k^*}\sqrt{2\varepsilon_k}\right\} + \frac{9}{n^2 k^2}$$

[since $\mathbb{P}(A_k \cap B_k)^c \leq 9/(n^2 k^2)$ by (4.6) in Lemma 4.3]

$$\leq \mathbb{P}\left\{\sup_{f \in \mathcal{F}_k^*} |L(f) - \widehat{L}(f)| \geq 4\mathbb{E}\widehat{R}_{\mathcal{F}_k^*} + 2\varepsilon_k + \Sigma_{\mathcal{F}_k^*}\sqrt{2\varepsilon_k}\right\} + \frac{10}{n^2 k^2}$$

(by Proposition 4.4)

$$\leq \mathbb{P}\left\{\sup_{f \in \mathcal{F}_k^*} |L(f) - \widehat{L}(f)| \geq 2\mathbb{E}\sup_{f \in \mathcal{F}_k^*} |\widehat{L}(f) - L(f)| + \frac{4\varepsilon_k}{3} + \Sigma_{\mathcal{F}_k^*}\sqrt{2\varepsilon_k}\right\}$$
$$+ \frac{10}{n^2 k^2} \qquad [\text{by (4.7)}]$$
$$\leq \frac{11}{n^2 k^2} \qquad \text{(by Proposition 4.5)}.$$

This inequality and Lemma 2.1 imply the first assertion of the theorem. The other statements–the probability bound and the almost-sure statement– follow by invoking Lemma 2.2 and the preceding argument, which also shows that

$$\mathbb{P}\{\widehat{C}_k \leq (L - \widehat{L})(f_k^*)\} \leq \frac{11}{n^2 k^2},$$

although the last assertion could be shown in a much easier way as it only involves a single function $f_k^*$. The proof of Theorem 4.1 is complete. □

In the proof of Theorem 4.2 we need the symmetrization device

$$(4.8) \qquad \mathbb{E}\widehat{R}_{\mathcal{F}_k} \leq 2\mathbb{E}\sup_{f \in \mathcal{F}_k} |\widehat{L}(f) - L(f)| + \frac{\sup_{f \in \mathcal{F}_k} L(f)}{\sqrt{n}}$$

(see, e.g., [20], page 18), and also the following result due to Massart [18]. (The version stated here is taken from [16].)

PROPOSITION 4.6. *Set* $\Sigma_k = \sup_{f \in \mathcal{F}_k} \sqrt{L(f)(1 - L(f))}$. *Then, for all* $n \geq 1$,

$$\mathbb{E}\sup_{f \in \mathcal{F}_k} |\widehat{L}(f) - L(f)| \leq \frac{8\mathbb{E}\log 2\mathbb{S}_k(X_1^{2n})}{n} + 4\sqrt{\frac{2\Sigma_k^2 \, \mathbb{E}\log 2\mathbb{S}_k(X_1^{2n})}{n}}.$$

PROOF. The statement follows almost immediately from Theorem 1.10 in [16] by noting that the worst-case shatter coefficients may be replaced with impunity by the random shatter coefficients. □

PROOF OF THEOREM 4.2. Observe that on the event $A_k \cap B_k$, $\widehat{\mathcal{F}}_k \subseteq \overline{\mathcal{F}}_k$, where $\overline{\mathcal{F}}_k$ is as defined in Theorem 4.2. Indeed, for any $f \in \widehat{\mathcal{F}}_k$,

$$\begin{aligned}
L(f) &\leq 2\widehat{L}(f) + u_k && \text{(by definition of } A_k\text{)} \\
&\leq 2[16\widehat{L}(\hat{f}_k) + 15\hat{u}_k] + u_k && \text{(by definition of } \widehat{\mathcal{F}}_k\text{)} \\
&\leq 32\widehat{L}(\hat{f}_k) + 31\bar{u}_k && \text{(by definition of } B_k\text{)} \\
&\leq 32\widehat{L}(f_k^*) + 31\bar{u}_k && \text{(by definition of } \hat{f}_k\text{)} \\
&\leq 32[2L_k^* + u_k] + 31\bar{u}_k && \text{(by definition of } A_k\text{)} \\
&= 64L_k^* + 63\bar{u}_k.
\end{aligned}$$

Also, we notice that on the event $A_k$,

$$\widehat{L}(\hat{f}_k) \leq \widehat{L}(f_k^*) \leq 2L_k^* + u_k.$$

These observations imply that

$$\begin{aligned}
\widehat{C}_k I_{A_k \cap B_k} &\leq 8\widehat{R}_{\overline{\mathcal{F}}_k} + 10\varepsilon_k + 2\sqrt{64L_k^* + 63\bar{u}_k}\sqrt{2\varepsilon_k} \\
&\leq 8\widehat{R}_{\overline{\mathcal{F}}_k} + 10\varepsilon_k + 16\sqrt{L_k^* + \bar{u}_k}\sqrt{2\varepsilon_k}.
\end{aligned}$$

Consequently, it follows from Lemma 4.3 that

$$\begin{aligned}
\mathbb{E}\widehat{C}_k &\leq \mathbb{E}\widehat{C}_k I_{A_k} + \mathbb{P}(A_k \cap B_k)^c \\
&\leq 8\mathbb{E}\widehat{R}_{\overline{\mathcal{F}}_k} + 10\varepsilon_k + 16\sqrt{L_k^* + \bar{u}_k}\sqrt{2\varepsilon_k} + 9(nk)^{-2} \\
&\leq 8\mathbb{E}\widehat{R}_{\overline{\mathcal{F}}_k} + 15\varepsilon_k + 16\sqrt{L_k^* + \bar{u}_k}\sqrt{2\varepsilon_k}.
\end{aligned}$$

This bound and Theorem 4.1 yield the first inequality of Theorem 4.2. The second inequality follows from the symmetrization (4.8) and Proposition 4.6. □

**Acknowledgments.** We thank Olivier Bousquet for his invaluable remarks and advice. We also appreciate the helpful remarks of two referees.

DEPARTMENT OF ECONOMICS
POMPEU FABRA UNIVERSITY
RAMON TRIAS FARGAS 25-27
08005 BARCELONA
SPAIN
E-MAIL: lugosi@upf.es

DEPARTMENT OF STATISTICS
FLORIDA STATE UNIVERSITY
TALLAHASSEE, FLORIDA 32306-4330
USA
E-MAIL: wegkamp@stat.fsu.edu